\documentclass[10pt,oneside]{jm}
\usepackage{amsfonts,amssymb,amsmath,rotating}
\usepackage[backref, pageanchor=false, colorlinks=true,%
            linkcolor=blue, citecolor=red, urlcolor=magenta]{hyperref}

\newcommand{\abs}{\vskip 0.5em\noindent\rm}
\newcommand{\Abs}{\paragraph{}\hspace{-1em}\rm}
\newcommand{\AbsT}[1]{\paragraph{\hspace{-1em} #1}\rm}

\newcommand{\bfi}{\noindent {\bf i)} }

\newcommand{\bfii}{\noindent {\bf ii)} }

\newcommand{\bfiii}{\noindent {\bf iii)} }

\newcommand{\bfa}{\noindent {\bf a)} }
\newcommand{\bfb}{\noindent {\bf b)} }
\newcommand{\bfc}{\noindent {\bf c)} }
\newcommand{\bfd}{\noindent {\bf d)} }

\newcommand{\absr}{\abs\hrulefill}

\newcommand{\C}{\mathbb C}

\newcommand{\F}{\mathbb F}

\newcommand{\K}{\mathbb K}
\newcommand{\N}{\mathbb N}

\newcommand{\Q}{\mathbb Q}

\newcommand{\Z}{\mathbb Z}

\newcommand{\cI}{\mathcal I}

\newcommand{\cO}{\mathcal O}

\newcommand{\cQ}{\mathcal Q}
\newcommand{\ccR}{\mathcal R}

\newcommand{\dt}{\delta}

\newcommand{\ph}{\varphi}

\newcommand{\End}{\text{End}}

\newcommand{\Irr}{\text{Irr}}

\newcommand{\rad}{\text{rad}}

\newcommand{\Stab}{\text{Stab}}

\newcommand{\Tr}{\text{Tr}}

\newcommand{\baraut}{\ov{\text{\hspace{0.5em}\rule{0em}{0.55em}}}}

\newcommand{\ld}{,\ldots\hskip0em ,}
\newcommand{\lr}[1]{\langle #1\rangle}

\newcommand{\ov}[1]{\overline{#1}}

\newcommand{\cl}[1]{\lceil #1\rceil}

\newcommand{\wh}[1]{\widehat{#1}}
\newcommand{\mt}{\mapsto}

\newcommand{\ra}{\rightarrow}

\newcommand{\llra}{\longleftrightarrow}

\newcommand{\cn}{\colon}
\newcommand{\dcup}{\stackrel{.}{\cup}}

\newcommand{\spmid}{\,\mid\,}
\newcommand{\spnmid}{\,\nmid\,}

\newcommand{\sseq}{\subseteq}

\newcommand{\tm}{\times}
\newcommand{\otm}{\otimes}

\newcommand{\MAGMA}{{\sf MAGMA}}
\newcommand{\GAP}{{\sf GAP}}

\newcommand{\MA}{{\sf MeatAxe}}
\newcommand{\IMA}{{\sf IntegralMeatAxe}}
\newcommand{\ORB}{{\sf ORB}}

\begin{document}
\raggedbottom
\pagestyle{myheadings}
\markboth{}{}
\thispagestyle{empty}
\setcounter{page}{1}

\begin{center} \Large\bf
Is the projective cover of the trivial module \\
in characteristic $11$ for the sporadic simple \\
Janko group $J_4$ a permutation module? \vspace*{1em} \\
\large\rm Jürgen Müller \vspace*{1em} \\
\large\it Dedicated to the memory of Richard Parker. \vspace*{1.5em} \\
\end{center}

\begin{abstract} \vspace*{-3em}\noindent\hrulefill\vspace*{2em} \\ \noindent 
We determine the ordinary character of the projective cover of the trivial 
module in characteristic $11$ for the sporadic simple Janko group $J_4$,
and answer the question posed in the title. 

\vspace*{0.5em} \noindent

\vspace*{0.5em} \noindent 

\vspace*{0em} \noindent \hrulefill 
\end{abstract}


\section{Introduction}

\Abs
Let $p$ be a prime, let $\F$ be a field of characteristic $p$,
let $G$ be a finite group, and let $P_{\F_G}$ be the projective 
cover of the trivial $\F[G]$-module $\F_G$.

\abs
The present article is motivated by the recent paper \cite{MalRob},
which deals with the question when $P_{\F_G}$ is a permutation module. 
This amounts to asking whether there is a subgroup $H\leq G$ such that 
$P_{\F_G}$ is isomorphic to the induced module $\F_H^G$.
Obviously, this is the case whenever $G$ is a $p'$-group, so that we 
may assume that $p\spmid |G|$. Moreover, by \cite[Cor.2.6]{MalRob}, 
if $G$ has the above property, then so has any composition factor of $G$,
which shifts focus to non-abelian simple groups. 

\abs 
Now, in \cite{MalRob}, the non-abelian simple groups having the above property 
are classified, apart from the groups of Lie type in defining characteristic 
and, amongst the sporadic simple groups, the largest Janko group $J_4$ 
in characteristic $11$.

\abs
The latter escapes all purely character theoretic attacks. Fortunately, 
I have been able to provide computational help to answer the above question: 
It turns out that $P_{\F_{J_4}}$ is {\it not} a permutation module, where
$\F$ is a field of characteristic $11$. This is already reported in 
\cite[Thm.4.1]{MalRob}, albeit without proof. It is the purpose of the 
present article to provide the details of the computations I made. 

\abs
Actually, these computations even reveal the projective indecomposable
(ordinary) character $\Psi_{\F_{J_4}}$ afforded by $P_{\F_{J_4}}$, 
and a few more projective indecomposable characters belonging to the 
principal $11$-block of $J_4$. About the latter virtually nothing is known 
so far, so that the results presented here might be the first steps towards
finding its decomposition matrix.

\Abs
We describe the general approach, see also \cite[Ch.2]{MalRob}:
Note first that, given $G$, the above property only depends on $p$, but
not on the particular choice of $\F$, so that later on
we will let $\F$ be the prime field $\F_p$.

\abs
Now, if $G$ has the above property, since $P_{\F_G}$ is indecomposable,
it necessarily is a transitive permutation module. Thus it has shape
$\F_H^G$, for some $H\leq G$, where since $P_{\F_G}$ is projective 
$H$ necessarily is a $p'$-subgroup. Moreover, since for any $p'$-subgroup
$U\leq G$ the module $\F_U^G$ has $P_{\F_G}$ as a direct summand, 
we conclude that $H$ has maximal order amongst all $p'$-subgroups of $G$.

\abs
If the $p$-modular decomposition matrix of $G$ is known,
then by the theory of trivial-source modules, see \cite[Ch.II.12]{Landrock},
we may just check whether the permutation character $1_H^G$ 
associated with the action of $G$ on the cosets of $H$
coincides with the projective indecomposable character $\Psi_{\F_G}$. 
But, as mentioned above, for the case of interest 
to us we are not at all in this comfortable position.

\abs
Hence we have to proceed otherwise: 
We pick any $p'$-subgroup $H\leq G$ of maximal order, 
and check whether $\F_H^G$ is an indecomposable $\F[G]$-module, where 
the latter property is equivalent to the endomorphism algebra 
$\End_{\F[G]}(\F_H^G)$ being a local $\F$-algebra. 
Thus, to pursue this, we have to analyze the structure of the 
endomorphism algebra of a permutation module, which for the group
in question due to sheer size is not too easy to handle.

\Abs
The present article is organized as follows:
In Section \ref{sec1} we describe the background concerning
endomorphism algebras and the orbit enumeration techniques used; 
in Section \ref{sec2} we provide some character theoretic data on the 
principal $11$-block of $G:=J_4$, and specify the subgroup $H$ to consider;
in Section \ref{sec3} we consider the $G$-action on the set $\cO$ of 
cosets of $H$, and apply \ORB{} to find the $H$-orbits in $\cO$; 
in Section \ref{sec4} we consider the endomorphism algebra $E$ of the 
permutation module afforded by $\cO$, and determine the character table 
of $E$;
and finally in Section \ref{sec5} we compute the decomposition matrix of 
$E$, and answer the question we started with. 

\AbsT{Acknowledgments.}
It is a great honor to have this opportunity to thank Richard Parker 
for a wealth of mathematical ideas he was always keen to share 
with everybody. In particular, the present article owes much
to his work, as a glance into the list of references reveals. 
Notably, proving the sheer existence of $J_4$ was Richard Parker's 
original motivation to invent the \MA{}.


\abs
Moreover, I would also like to thank Thomas Breuer for inspiring discussions
(not only) about the topic of this article, and the referee for their
careful reading.

\section{Prerequisites}\label{sec1}

\AbsT{Endomorphism algebras.}
We recall the necessary facts about the structure of endomorphism algebras 
of permutation modules, thereby fixing the notation used later; as a 
general reference see \cite[Ch.II.12]{Landrock}, while the background of the 
particular approach we follow is described in detail in \cite{Habil}:

\abs\bfa
Let $G$ be a finite group, let $H\leq G$ be a subgroup, let $\cO$ 
be the set of (right) cosets of $H$ in $G$, and let $n:=|\cO|$.
Moreover, if $R$ is an integral domain, let $R_H^G$ be the 
permutation $R[G]$-module associated with $\cO$, and let 
$E_R:=\End_{R[G]}(R_H^G)$ be its endomorphism $R$-algebra.
Then $E_R$ is $R$-free of rank $r:=\lr{1_H^G,1_H^G}_G$, where 
$1_H^G$ is the induced character from the trivial character $1_H$ of $H$,
and $\lr{\cdot,\cdot}_G$ denotes the usual scalar product 
on the characters of $G$.

\abs
Let $\{v_1\ld v_r\}$ be a set of representatives of the $H$-orbits 
$\cO_j:=(v_j)^H\sseq\cO$, where $v_1$ denotes the coset $H$ itself, and
let $g_j\in G$ such that $v_1g_j=v_j$. Moreover, let $H_j:=\Stab_H(v_j)$, 
and let $n_j:=|\cO_j|$. The {\it paired orbit} $\cO_{j^\ast}$ of $\cO_j$ 
is defined to be the $H$-orbit containing $v_1g_j^{-1}$.

\abs
Let $\cO_j^+:=\sum_{v\in\cO_j}v\in R_H^G$ be the associated {\it orbit sum}. 
Then $E_R$ has a distinguished $R$-basis $\{A_1\ld A_r\}$, being called its
{\it Schur basis}, where $A_j$ is defined by $v_1\mt\cO_j^+$, and
subsequent extension by $G$-transitivity to all of $\cO$;
in particular $A_1$ is the identity map. 
Thus, abbreviating $E=E_\Z$, we have $E_R=E\otm_\Z R$.

\abs
Writing $A_iA_j=\sum_{k=1}^r p_{ijk}A_k\in E$, 
the associated (non-negative) structure constants, also being called 
{\it intersection numbers}, are given as
$p_{ijk}=\frac{n_i}{n_k}\cdot c_{jk}(g_i)\in\Z$,
using the {\it orbit counting numbers} 
$c_{jk}(g_i):=|\cO_jg_i\cap\cO_k|\in\N_0$.
Thus the (right) regular representation of $E$, with respect to its
Schur basis, is given by the {\it intersection matrices}
$P_j:=[p_{ijk}]_{ik}\in\Z^{r\tm r}$. In particular, the first row
of $P_j$ is given by $p_{1jk}=\dt_{jk}$, that is, consists of the $j$-th
unit vector.

\abs\bfb
Let $K$ be a field. Then, for any $E_K$-module $V$,
the trace map $\ph_V\cn E_K\ra K\cn A\mt\Tr_V(A)$
is called the {\it character} afforded by $V$. Letting
$\Irr(E_K):=\{\ph_1\ld\ph_s\}$ be the set of characters afforded by the 
irreducible $E_K$-modules, we obtain the 
{\it character table} $\Phi:=[\ph_i(A_j)]_{ij}\in K^{s\tm r}$.

\abs
Since $\C_H^G$ is a semi-simple $\C[G]$-module, 
$E_\C$ is a (split) semi-simple $\C$-algebra, and we have a natural
bijection between the irreducible representations of $E_\C$ and the 
distinct constituents of $\C_H^G$, being called {\it Fitting correspondence};
in terms of irreducible characters the Fitting correspondent of 
$\ph\in\Irr(E_\C)$ is denoted by $\chi_\ph$.
Moreover, $E_\C$ is commutative if and only if $\C_H^G$ is multiplicity-free.

\abs
We have $\ph(A_1)=m_{\chi_\ph}=\lr{1_H^G,\chi_\ph}_G$, 
the multiplicity of $\chi_\ph$ as a constituent of $1_H^G$. 
The Fitting correspondent $\ph_1\in \Irr(E_\C)$ of $1_G$ is given by 
$\ph_1(A_j)=n_j$; it is the only irreducible
character of $E_\C$ whose values on the Schur basis consist of non-negative 
integers only. We have the following {\it orthogonality relations}
between characters $\ph,\ph'\in\Irr(E_\C)$, where $\baraut$ denotes complex 
conjugation, and where we have $\ov{\ph(A_j)}=\ph(A_{j^\ast})$:
$$ \frac{1}{n}\cdot\sum_{j=1}^r\frac{1}{n_j}\cdot\ov{\ph(A_j)}\cdot\ph'(A_j)
=\dt_{\ph,\ph'}\cdot \frac{m_{\chi_\ph}}{\chi_\ph(1)} .$$

\abs\bfc
The endomorphism algebra $E$ admits a decomposition theory,
similar to the one for group algebras: 
Let $\ccR$ be a discrete valuation ring in an algebraic number field $\K$,
such that the maximal ideal $\wp\lhd\ccR$ contains $p$, and let 
$\F:=\ccR/\wp$. In practice, in order to keep data consistent, 
we make the same conventional choices for $\ccR$ and $\wp$ as in 
\cite{ModAtlas}.
Moreover, we assume that $\K$ and $\F$ are large enough so that both  
$E_\K$ and $E_\F$ are split.

\abs
Then any finitely generated $E_\K$-module can be realized by an
$E_{\ccR}$-lattice $V$, and {\it $\wp$-modular reduction},
mapping $V$ to $V_\F:=V\otm_{\ccR}\F$, yields
a $\Z$-linear decomposition map $D_\wp\cn G(E_\K)\ra G(E_\F)$
between the associated Grothendieck groups. Its matrix
with respect to the $\Z$-bases consisting of the respective irreducible
representations is called the associated {\it decomposition matrix}.

\abs
Since $\K$ is a splitting field for $E_\K$, we have $\Irr(E_\C)=\Irr(E_\K)$,
which is $\K$-linearly independent, so that we may identify $G(E_\K)$ with 
$\Z\Irr(E_\C)$. Since $\F$ is a splitting field for $E_\F$, similarly 
$\Irr(E_\F)$ is $\F$-linearly independent. 
Since $\Irr(E_\K)$ has values in $\ccR$, for any element of $E_\ccR$,
we conclude that $\wp$-modular reduction induces a $\Z$-linear map 
$D_\wp\cn\Z\Irr(E_\C)\ra\F\Irr(E_\F)$.

\abs\bfd
Let $S$ be a simple $E_\F$-module, with associated projective 
indecomposable module $P_S\cong e_S E_\F$, for some suitable idempotent 
$e_S\in E_\F$. Moreover, for $\ph\in\Irr(E_\K)$ let $V_\ph$ be an  
$E_{\ccR}$-lattice such that $(V_\ph)_\K:=V_\ph\otm_\ccR\K$ has 
character $\ph$, and let $e_\ph\in E_\K$ be an idempotent such that 
$e_\ph E_\K\cong (V_\ph)_\K$.

\abs
Any idempotent $e\in E_\F$ can be {\it lifted} to $E_\ccR$, that is, 
there is an idempotent $\wh e\in E_\ccR$ such that $\wh e\otm 1_\F=e\in E_\F$.
In particular, there is a projective indecomposable $E_\ccR$-lattice 
$\wh P_S\cong\wh e_S E_\ccR$ lifting $P_S$. Thus for the multiplicity 
of $S$ as a constituent of $V_\F$ we have {\it Brauer reciprocity} 
$[V_\F\cn S]=[(\wh P_S)_\K\cn V_\K]$, and for the 
{\it Cartan numbers} of $E_\F$ we have $[P_S\cn S']=[P_{S'}\cn S]$.
In particular, since $E_\F\cong\bigoplus_S (P_S)^{\oplus\dim_\F(S)}$ 
as $E_\F$-modules, this entails
$$ \dim_\F(P_S)=\sum_{S'}\dim_\F(S')\cdot[P_S\cn S']
=\sum_{S'}\dim_\F(S')\cdot[P_{S'}\cn S]=[E_\F\cn S] .$$

\abs
This relates to Fitting correspondence as follows: For an irreducible 
character $\chi$ of $G$ occurring as a constituent of $1_H^G$, let 
$V_\chi$ be an $\ccR[G]$-lattice such that $(V_\chi)_\K$ has character $\chi$.
Then we have $(\K_H^G)e_\ph\cong(V_{\chi_\ph})_\K$ as $\K[G]$-modules,
and thus 
$$ [(\K_H^G)\wh e_S\cn(V_{\chi_\ph})_\K] 
=[\wh e_S E_\K\cn(V_\ph)_\K]=[(\wh P_S)_\K\cn(V_\ph)_\K]=[(V_\ph)_\F\cn S] .$$

\AbsT{Enumeration of long orbits.}
To facilitate computations with (large) permutation representations 
we use the \GAP{} package \ORB{} \cite{ORB}, where its orbit enumeration
techniques are described comprehensively in \cite{Habil}, and an extended
worked application is presented in \cite{BMpaper}. 
We give a brief sketch of the approach:

\abs
Let $G$ be a (large) finite group, and let $\cO$ be a (large) 
transitive $G$-set, which we assume to be implicitly given,
for example as a $G$-orbit of a vector $v_1$ in an $F[G]$-module $V$ over 
a finite field $F$. Letting $H\leq G$ be a (still large) subgroup,
we are interested in classifying the $H$-orbits $\cO_j$ in $\cO$, finding 
their length $n_j$, representatives $v_j\in\cO_j$, elements $g_j\in G$ 
such that $v_1\cdot g_j=v_j$, and the stabilizers $H_j=\Stab_H(v_j)$. 
To achieve this, we assume to be able to compute efficiently within $H$ 
(but not within $G$), for example by having a (smallish) faithful permutation
representation of $H$ at hand.

\abs
To find the $H$-orbits in $\cO$, we choose a (smallish) helper subgroup 
$K\leq H$, and enumerate the various $H$-orbits $\cO_j\sseq\cO$ by the
$K$-orbits they contain. To do so, we choose a (not too small)
helper $K$-set $\cQ$ together with a homomorphism $\pi_K\cn\cO\ra\cQ$ 
of $K$-sets, which again we assume to be implicitly given, for example
by an $F[K]$-quotient module of $V$.

\abs
Moreover, we assume that $K$ has sufficiently long orbits in $\cQ$, 
and that we are able to classify them, by giving representatives,
their stabilizers in $K$, as well as complete Schreier trees.
(Thus for the $K$-action on $\cQ$ we are facing a similar problem 
as for the $H$-action on $\cO$, apart from the requirement 
on Schreier trees. So we could just recurse. Actually, the full 
functionality of \ORB{} supports this, but for the present purposes 
we will get away with a single helper subgroup.)

\abs
For any $K$-orbit in $\cQ$, the chosen representative is called its
{\it distinguished point} (although it might be chosen arbitrarily).
Then, for any $K$-orbit $\cO'\sseq\cO$, the $\pi_K$-preimages of the 
distinguished point of $\pi_K(\cO')\sseq\cQ$ are likewise called the 
distinguished points of $\cO'$.
Hence, to enumerate an $H$-orbit $\cO_j$ by enumerating the $K$-orbits 
contained in it, we only have to store the associated distinguished points, 
and a Schreier tree telling us how to reach them from $v_j$.

\abs 
For any $H$-orbit $\cO_j$ we are content of finding only as many 
$K$-orbits contained in it as are needed to cover more than half of it;
this is equivalent to knowing $n_j$ and $|H_j|$. Then we have a 
randomized membership test for $\cO_j$, and a deterministic test to 
decide whether the $H$-orbits found are actually pairwise disjoint.

\abs
The number of points of $\cO_j$ covered by the above enumeration process,
divided by the number of distinguished points actually stored is called 
the {\it saving factor} achieved. The maximum saving factor possible is 
$|K|$, which is achieved if and only if $K$ has only regular orbits in the
$\pi_K$-image of the part of $\cO_j$ covered.

\AbsT{Computational tools.}
To facilitate group theoretic and character theoretic computations we use 
the computer algebra system \GAP{} \cite{GAP}, its comprehensive database 
{\sf CTblLib} \cite{CTblLib} of ordinary and modular character tables, 
and its library {\sf TomLib} \cite{TomLib} of tables of marks.
In particular, {\sf CTblLib} encompasses the data given in the
{\sf Atlas} \cite{Atlas} and in the {\sf ModularAtlas} \cite{ModAtlas},
as well as the additional data collected on the 
{\sf ModularAtlasHomepage} \cite{ModAtlasH}.

\abs
As far as matrix representations over finite fields are concerned,
we use the \MA{} \cite{MA}, whose basic ideas go back to \cite{Parker}, 
where we also use its extensions to compute submodule 
lattices \cite{LMR} and direct sum decompositions \cite{LuxSz}.

\abs
Computations with matrix representations over the integers and over 
the rational numbers are facilitated by the \GAP{} package \IMA{} \cite{IMA}, 
which is developed and used heavily in \cite{GeMu}, but owes much to 
\cite{Parker2}. (The \IMA{} package is as yet unpublished, 
but I am of course happy to provide the code to everybody interested.
Moreover, as an alternative, similar functionality is available in the
computer algebra system \MAGMA{} \cite{MAGMA}.)

\abs
Data concerning explicit permutation representations,
ordinary and modular matrix representations, and the embedding of 
(maximal) subgroups of sporadic simple groups is available in the 
{\sf AtlasOfGroupRepresentations} \cite{AtlasGrpRep}, and through the 
\GAP{} package {\sf AtlasRep} \cite{AtlasRep}. 
For a wealth of group theoretical information used throughout we refer to 
the {\sf Atlas} \cite{Atlas}, whose notational conventions we follow; 
in particular we let $r_n:=\sqrt{n}$ be the positive square root of $n\in\N$.

\section{The principal $11$-block of $J_4$}\label{sec2}

\Abs
From now on let $G:=J_4$.

\abs
Then $G$ has the principal block $B_0$ as its only $11$-block of 
positive defect. The defect groups of $B_0$, that is, the Sylow 
$11$-subgroups of $G$, are extraspecial of shape $11^{1+2}_+$, and have 
the rare property of being trivial-intersection subgroups.
There are $k_0:=49$ irreducible ordinary characters and 
$l_0:=40$ irreducible modular characters belonging to $B_0$.

\abs
Sadly enough, this is virtually all what is known about the 
decomposition numbers of $B_0$, according to the {\sf ModularAtlasHomepage},
where $B_0$ is a prominent gap, in particular in view
of the trivial-intersection property of its defect groups.

\Abs
Using the conjugacy classes of maximal subgroups of $G$, as reproduced 
in the {\sf Atlas}, it turns out that the unique class of subgroups
of $11'$-subgroups of maximal order is given by the maximal subgroups 
of $G$ of shape $2^{10}\cn L_5(2)$. 
Let $H<G$ be a representative of this class.

\abs
We have $|H|=10.239.344.640$ and $[G\cn H]=8.474.719.242$.
The decomposition of the permutation character $1_H^G$ into the
irreducible ordinary characters $\chi_i$ of $G$ is given in 
Table \ref{pchtbl}, where the $\chi_i$ are ordered as in the {\sf Atlas},
we indicate generators of their (quadratic) character fields, and 
$m_i:=\lr{1_H^G,\chi_i}_G$. In particular, we have $r:=\lr{1_H^G,1_H^G}_G=27$,
and all constituents of $1_H^G$ belong to $B_0$. 

\abs
All constituents except $\chi_{19/20}$ are $11$-rational characters, 
being fixed by the Frobenius automorphism. Amongst them, $\chi_{23/24}$, 
$\chi_{36/37}$ and $\chi_{38/39}$ are pairs of Galois conjugate characters.
The constituents $\chi_{19/20}$ are non-$11$-rational, Galois conjugate 
characters, restricting to the same character on $11$-regular classes.

\begin{table}\caption{The permutation character $1_H^G$.}\label{pchtbl}
$$ \begin{array}{|r|c||c|} \hline 
\chi_i & \Q(\chi_i) & m_i \\ \hline \hline
 1 & & 1 \\
 8 & & 1 \\
11 & & 1 \\
14 & & 1 \\
19 & r_{33} & 2 \\
20 & r_{33} & 2 \\
\hline \end{array} \quad\quad\quad
\begin{array}{|r|c||c|} \hline 
\chi_i & \Q(\chi_i) & m_i \\ \hline \hline
21 & & 2 \\
22 & & 1 \\
23 & r_3 & 1 \\
24 & r_3 & 1 \\
29 & & 1 \\
30 & & 1 \\
\hline \end{array} \quad\quad\quad
\begin{array}{|r|c||c|} \hline 
\chi_i & \Q(\chi_i) & m_i \\ \hline \hline
32 & & 1 \\
36 & r_5 & 1 \\
37 & r_5 & 1 \\
38 & r_5 & 1 \\
39 & r_5 & 1 \\
51 & & 1 \\
\hline \end{array} $$
\absr\end{table}

\Abs\label{psi}
The projective indecomposable character $\Psi_{\F_G}$ is a summand of 
$1_H^G$. Thus, writing $\Psi_{\F_G}=\sum_{i\in\cI} d_i\chi_i$, where
$d_i\in\N_0$ and $\cI$ is the index set occurring in the first column
of Table \ref{pchtbl}, we have $0\leq d_i\leq m_i$;
in particular $d_1=m_1=1$.

\abs
Moreover, we have $\Psi_{\F_G}(g)=0$ whenever $g\in G$ is $11$-singular.
Enforcing these conditions, it turns out that from the 
$2^{15}\cdot 3^3=884736$ candidates for $\Psi_{\F_G}$ allowed by the 
above inequalities, there are just $75$ admissible candidates left.

\abs
The above conditions say that $\Psi_{\F_G}$ is a generalized
projective character. Thus, by \cite[Cor.2.17, La.2.21]{Navarro},
we conclude that $\frac{1}{|C_G(g)|_{11}}\cdot\Psi_{\F_G}(g)$ 
is an algebraic integer, for any $g\in G$, in particular entailing that 
$11^3\spmid\Psi_{\F_G}(1)$.
Moreover, it turns out that they already imply $d_{19}=d_{20}$,
and neither complex conjugation nor the Frobenius automorphism 
yield further conditions. (Note that $\Psi_{\F_G}$ is not necessarily
a rational character, although the trivial character $1_G$ is.)


\abs
Unfortunately, we have not been able to find further purely
character-theoretic conditions to narrow down the set of 
admissible candidates for $\Psi_{\F_G}$.
In particular, neither restriction to maximal subgroups of $G$
and decomposition into their projective indecomposable characters 
(providing lower bounds on the $d_i$), nor induction of the 
projective cover of the trivial module of maximal subgroups
(providing upper bounds on the $d_i$) did yield any improvement.
At this point we have decided to revert to explicit computations, 
in particular applying \ORB{}.

\section{The permutation action}\label{sec3}

\Abs
To do explicit computations, we pick the $112$-dimensional (absolutely)
irreducible representation of $G$ over $\F_2$ from the {\sf AtlasRep} 
database. The latter is given in terms of (two) standard generators, 
in the sense of \cite{Wilson}. Words in the standard generators providing 
generators of a maximal subgroup $2^{10}\cn L_5(2)\cong H<G$ are also 
available in the {\sf AtlasRep} database.

\abs
Let $V\cong\F_2^{112}$ be the module underlying the above representation
of $G$. Using the \MA{}, it turns out that $H$ possesses a $1$-dimensional 
fixed space in $V$. Hence letting $v_1\in V$ be the unique non-zero
$H$-fixed vector, the $G$-action on the orbit $\cO:=(v_1)^G\sseq V$ 
is equivalent to its action on the cosets of $H$. Hence this provides 
an implicit realization of the latter action.

\abs
To store a vector in $V$, including header information, we need 
$(\cl{\frac{112}{8}}+4)\text{ Byte}=18$ Byte. Thus to store $\cO$ completely
we would need at least
$$ ([G\cn H]\cdot 18)\text{ Byte}=152.544.946.356\text{ Byte}
   \sim 150 \text{ GB} $$ 
of memory space, plus some more header information.
Although this would in principle be feasible nowadays, 
in view of the computations we are going to make within $\cO$,
we apply \ORB{}, trying to achieve a saving factor of $\sim 150$, say.
We set up the required framework:

\Abs\bfi
A faithful permutation representation of $H$ is found as follows:
Using the \MA{}, we determine the submodule lattice of the 
restriction $V_H$ of $V$ to $H$. It turns out that it possesses
a unique $16$-dimensional $\F_2[H]$-submodule $U$.  Moreover,
there is a faithful $H$-orbit of vectors in $U$ of length $310$. 
Computing the $H$-action on the latter yields an explicit permutation 
representation of $H$. 

\abs\bfii
Before actually choosing a helper subgroup, we already look for a helper 
$H$-set, being an epimorphic image of $\cO$ as $H$-sets:

\abs 
Since $V$ is a self-contragredient $\F_2[G]$-module, $V_H$ has a unique 
$16$-dimensional quotient $\F_2[H]$-module $W$, being the dual of the 
submodule $U$ considered above, and thus in particular being a faithful
$\F_2[H]$-module as well. Having cardinality $2^{16}=65536$, it is 
small enough to enumerate all its elements explicitly.

\abs
Using the \MA{}, we compute the natural quotient map $V_H\ra W$ of 
$\F_2[H]$-modules, which gives rise to a homomorphism $\pi_H\cn\cO\ra W$ 
of $H$-sets.  

\abs\bfiii
As helper subgroup we now choose $K:=N_H(31)\cong 31\cn 5$, 
the normalizer in $H$ of a Sylow $31$-subgroup, having order $|K|=155$. 
Words in the chosen generators of $H$ providing generators of $K$ are 
found with the help of \GAP{}, employing the permutation representation
of $H$ constructed above.

\abs
Returning to the matrix representation of $K$ on $V$, and going over 
to the quotient $W$, it turns out that the $K$-orbits of vectors in $W$ 
have lengths $[1^2,31^{14},155^{420}]$, where the exponents denote 
multiplicities. Hence the expected value of the length of the $K$-orbit 
of a randomly chosen vector in $W$ is $\sim 154$, so that we indeed expect
a suitable saving factor as envisaged above. 

\Abs 
We are now prepared to run \ORB{}, in order to find the decomposition
$\cO=\coprod_{j=1}^{27}\cO_j$ into $H$-orbits; recall that there are 
$r=27$ orbits indeed: 

\abs
We let $\cO_1:=(v_1)^H=\{v_1\}$. Then we randomly choose 
elements $g\in G$, and check whether $v_1g\in\cO$ belongs to 
one of the $H$-orbits already found. If not, then we have found 
a new $H$-orbit, $\cO_j$ say, and let $g_j:=g$ and $v_j:=v_1g$. 
We also consider $v_1g_j^{-1}\in\cO$ 
in order to detect non-self-paired $H$-orbits.

\abs
Letting this run for a certain while (actually some 2 hours on a single 
3 GHz CPU), we have been able to find $24$ of the $27$ orbits, making up 
all of $\cO$ up to $27001$ vectors. Thus only a fraction of 
$\sim 3.2\cdot 10^{-6}$ of $\cO$ is missing at this stage, 
making it highly improbable to conclude simply by random search.

\abs
Hence we set out to find the missing  (small) three orbits by identifying
the associated (large) point stabilizers: Using the library {\sf TomLib}
we find all subgroup orders of $L_5(2)$, and allowing for factors $2^i$, 
where $i\in\{0\ld 10\}$, we get a set of numbers encompassing the subgroup 
orders of $H$. Checking all $3$-tuples thereof whose indices in $H$ add up to 
$27001$ leaves the following candidates for the missing three orbit lengths:
$$ [31,930,26040], \quad [465,496,26040], \quad [31,7440,19530] .$$

\begin{table}\caption{The $H$-orbits in $\cO$.}
\label{orbtbl} 
$$ \begin{array}{|r|r||r|r|r||r|} \hline
j & j^\ast & n_j & |H_j| & H_j & \\
\hline \hline
 1 &    & 1 & 10239344640 &         2^{10}.L_5(2) &   1 \\
 2 &    & 31 &  330301440 &     2^{10}.2^4.L_4(2) &   1 \\
 3 &    & 930 &  11010048 & 2^9.2^6.(L_3(2)\tm 2) &   1 \\
 4 &    & 17360 &  589824 &         2^7.2^8.3^2.2 &   1 \\
 5 &    & 26040 &  393216 & 2^{1+12}.(D_8\tm S_3) &   1 \\
 6 &    & 27776 &  368640 &               2^9.S_6 &   1 \\
 7 &  8 & 416640 &  24576 &        2^5.2^6.D_{12} &   1 \\
 8 &  7 & 416640 &  24576 &        2^5.2^6.D_{12} &   1 \\
 9 &    & 624960 &  16384 &    (4\tm 2^4).2^4.2^4 &   1 \\
10 &    & 333120 &   3072 &       2^3.2^3.2^3.S_3 & 155 \\
11 &    & 4999680 &  2048 &      (D_8\tm 2^4).2^4 &   1 \\
12 & 13 & 6666240 &  1536 &             2^6.3.D_8 &   1 \\
13 & 12 & 6666240 &  1536 &             2^6.3.D_8 & 155 \\
14 &    & 9999360 &  1024 &    (D_8\tm D_8).2.2^3 & 155 \\
15 &    & 13332480 &  768 &               2^7.S_3 & 155 \\
16 &    & 53329920 &  192 &        2^2.2^2.D_{12} & 155 \\
17 &    & 66060288 &  155 &                  31.5 & 102 \\
18 & 19 & 79994880 &  128 &               2^3.2^4 & 155 \\
19 & 18 & 79994880 &  128 &               2^3.2^4 & 155 \\
20 &    & 159989760 &  64 &               2^3.2^3 & 155 \\
21 &    & 159989760 &  64 &               2^3.2^3 & 155 \\
22 &    & 319979520 &  32 &            D_8\tm 2^2 & 155 \\
23 &    & 341311488 &  30 &           D_{10}\tm 3 & 152 \\
24 &    & 1279918080 &  8 &                   2^3 & 155 \\
25 &    & 1279918080 &  8 &                   D_8 & 152 \\
26 &    & 2047868928 &  5 &                     5 & 150 \\
27 &    & 2559836160 &  4 &                     4 & 152 \\
\hline 
\end{array} $$ 
\absr\end{table}

\Abs
We proceed to decide which case occurs:
Recalling that $H\cong O_2(H)\cn L$, where $L\cong L_5(2)$, is a split 
extension, we have an embedding of $L$ into $H$. Since representatives 
of the conjugacy classes of subgroups of $L$ can be straightforwardly 
determined, this allows us to find representatives $S$ of the conjugacy 
classes of subgroups of $H$ of a fixed index. 

\abs
In turn, using the \MA{}, we compute the fixed space of $S$ on $V$, 
and check whether it contains a vector $v$ having an $H$-orbit of desired 
length and belonging to $\cO$. The latter property is verified by picking
random elements $g\in G$, and checking whether $vg\in V$ is contained in 
one of the known $H$-orbits, $\cO_j$ say. If so, then we may choose $v$ 
as a new $H$-orbit representative, and since we have enumerated $\cO_j$ 
already, we are readily able to find $h\in H$ such that $v_jh=vg$; 
thus we have $v_1\cdot g_jhg^{-1}=v$.

\abs\bfi
We first consider subgroups $S<H$ of index $31$. For these we have 
$O_2(H)<S$ and $S/O_2(H)\cong 2^4\cn L_4(2)$. There are 
two conjugacy classes. For one of them we are successful, excluding 
the second of the above cases.

\abs\bfii
Next, we consider subgroups $S<H$ of index $930$. Letting
$S^\ast:=S\cap O_2(H)$, we have either
$S^\ast=O_2(H)$ and $S/S^\ast\cong 2^6\cn L_3(2)$, or
$|S^\ast|=2^9$ and $S/S^\ast\cong 2^6\cn(L_3(2)\tm 2)$.
For one of the conjugacy classes of subgroups of the second shape 
we are successful, which brings us down to the first of the above cases.

\abs\bfiii
Finally, we consider subgroups $S<H$ of index $26040$, for which
we have $2^7\spmid |S^\ast|$. For one of the conjugacy classes of 
subgroups such that $|S^\ast|=2^7$ we are successful indeed. 
It turns out that in this case $S/S^\ast\cong 2^6\cn(D_8\tm S_3)$.

\abs
The results on the $H$-orbits in $\cO$ are collected in Table \ref{orbtbl}:
In the second column we indicate the non-self-paired orbits. 
In the fifth column we indicate the shape of $H_j$, where groups having
the same shape are actually isomorphic. (This is clear for 
$H_{7/8}$, $H_{12/13}$ and $H_{18/19}$ coming from paired orbits anyway.)
But, letting $H_j^\ast:=H_j\cap O_2(H)$, it turns out that $H_{7/8}^\ast$ 
and $H_{12/13}^\ast$ are different.

\abs
In the last column we also give the saving factors achieved for the
various $H$-orbits, which amounts to an average saving factor of $\sim 152$,
indicating that our choice of $K$ and $\cQ=W$ was not too bad. 
Still, we observe that the shorter $H$-orbits tend to have a saving factor 
of $1$, amounting to no saving at all. This could be due to our fairly 
ambitious choice of a helper $H$-set, rather than just a helper $K$-set, 
so that the quotient map might very well send a full $H$-orbit to the 
zero vector in $W$; but we have not analyzed this thoroughly.

\section{The endomorphism ring}\label{sec4}

\Abs
We consider the regular representation of the endomorphism algebra 
$E$ of $\Z_H^G$, which is $\Z$-free of rank $r=27$.
Computing the intersection matrices $P_j$ boils down to determining the 
orbit counting numbers $c_{jk}(g_i)=|\cO_jg_i\cap\cO_k|$: 

\abs
If $n_j$ is small enough so that $\cO_j$ can be enumerated explicitly, 
applying $g_i\in G$ and using the randomized orbit membership test in \ORB{},
we determine lower bounds for the numbers $c_{jk}(g_i)$. We are done once 
the figures we have found, running through all $k$, sum up to $n_j$. 
For example, $P_2$ is given in Table \ref{p2tbl}.

\abs
The orbits $\cO_j$ being ordered by increasing length, we successively 
compute $P_2,P_3,\ldots$, until the $\Q$-algebra generated by the 
matrices found so far has $\Q$-dimension $27$, and hence equals $E_\Q$.
Since $E$ is non-commutative, we need at least two generators, where it 
turns out that $P_3$ belongs to the $\Q$-algebra generated by $P_2$, but 
that $\{P_2,P_4\}$ suffice to generate $E_\Q$. 

\abs
The $P_j$ are associated with the regular representation of $E$,
with respect to its Schur basis. Hence the $\Q$-dimension of a 
candidate subalgebra as above can be determined by `spinning up' the 
first unit vector by applying the `standard basis algorithm', in the 
sense of \cite{Parker}, using the generators in question.
Moreover, since the first row of $P_j$ equals the $j$-th unit vector, 
decomposing it into the `standard basis' found above yields the complete 
intersection matrix $P_j$.  In practice, all of this is done using the \IMA{}.

\begin{table}\caption{The intersection matrix $P_2$.}\label{p2tbl}
\tiny \setlength{\arraycolsep}{0.44em}
$$ \left[\begin{array}{rrrrrrrrrrrrrrrrrrrrrrrrrrr}
.&1&.&.&.&.&.&.&.&.&.&.&.&.&.&.&.&.&.&.&.&.&.&.&.&.&. \\
31&.&1&.&.&.&.&.&.&.&.&.&.&.&.&.&.&.&.&.&.&.&.&.&.&.&. \\
.&30&2&.&1&.&.&.&.&.&.&.&.&.&.&.&.&.&.&.&.&.&.&.&.&.&. \\
.&.&.&1&4&.&.&1&.&.&.&.&.&.&.&.&.&.&.&.&.&.&.&.&.&.&. \\
.&.&28&6&2&.&.&.&1&.&.&.&.&.&.&.&.&.&.&.&.&.&.&.&.&.&. \\
.&.&.&.&.&1&.&2&.&.&.&.&.&.&.&.&.&.&.&.&.&.&.&.&.&.&. \\
.&.&.&.&.&.&1&.&4&3&.&.&.&.&.&.&.&.&.&.&.&.&.&.&.&.&. \\
.&.&.&24&.&30&.&.&.&.&1&1&.&.&.&.&.&.&.&.&.&.&.&.&.&.&. \\
.&.&.&.&24&.&6&.&2&.&1&.&.&1&.&.&.&.&.&.&.&.&.&.&.&.&. \\
.&.&.&.&.&.&24&.&.&.&.&.&.&.&1&.&.&.&1&.&.&.&.&.&.&.&. \\
.&.&.&.&.&.&.&12&8&.&5&.&.&.&3&.&.&1&.&.&.&.&.&.&.&.&. \\
.&.&.&.&.&.&.&16&.&.&.&6&.&.&.&.&.&.&.&1&.&.&.&.&.&.&. \\
.&.&.&.&.&.&.&.&.&.&.&.&1&4&.&.&.&2&.&.&.&.&.&.&.&.&. \\
.&.&.&.&.&.&.&.&16&.&.&.&6&2&.&.&.&.&1&.&1&.&.&.&.&.&. \\
.&.&.&.&.&.&.&.&.&4&8&.&.&.&3&.&.&.&.&.&.&1&.&.&.&.&. \\
.&.&.&.&.&.&.&.&.&.&.&.&.&.&.&1&.&4&.&.&.&.&.&1&.&.&. \\
.&.&.&.&.&.&.&.&.&.&.&.&.&.&.&.&.&.&.&.&.&.&.&.&.&1&. \\
.&.&.&.&.&.&.&.&.&.&16&.&24&.&.&6&.&.&.&1&1&1&.&1&.&.&. \\
.&.&.&.&.&.&.&.&.&24&.&.&.&8&.&.&.&.&5&.&.&2&.&.&1&.&. \\
.&.&.&.&.&.&.&.&.&.&.&24&.&.&.&.&.&2&.&1&4&4&.&.&2&.&. \\
.&.&.&.&.&.&.&.&.&.&.&.&.&16&.&.&.&2&.&4&1&.&.&1&.&.&1 \\
.&.&.&.&.&.&.&.&.&.&.&.&.&.&24&.&.&4&8&8&.&3&.&1&.&.&2 \\
.&.&.&.&.&.&.&.&.&.&.&.&.&.&.&.&.&.&.&.&.&.&1&4&4&.&. \\
.&.&.&.&.&.&.&.&.&.&.&.&.&.&.&24&.&16&.&.&8&4&15&3&4&5&4 \\
.&.&.&.&.&.&.&.&.&.&.&.&.&.&.&.&.&.&16&16&.&.&15&4&6&5&3 \\
.&.&.&.&.&.&.&.&.&.&.&.&.&.&.&.&31&.&.&.&.&.&.&8&8&10&8 \\
.&.&.&.&.&.&.&.&.&.&.&.&.&.&.&.&.&.&.&.&16&16&.&8&6&10&13 \\
\end{array}\right] $$
\absr\end{table}

\Abs
We proceed to determine the character table of $E_\C$ and the Fitting 
correspondence: To do so, we decompose the regular representation of
$E_\Q$, which is semi-simple, but not split. The ordinary constituents 
of $\Q_H^G$, as compared to those of $\C_H^G$ given in Table \ref{pchtbl}, 
tell us that $E_\Q$ has $14$ irreducible representations, of degree 
$[1^9,2^4,4] = [1^9,(1+1)^3,2,(2+2)]$, where exponents denote their 
multiplicity, and brackets denote their splitting over $\C$.

\abs
Using the \IMA{} and \GAP{}, we compute the characteristic polynomials
of $P_2$ and $P_4$, their factorization over $\Q$, and for the irreducible
divisors $f$ and $g$ occurring, respectively, we determine the 
$E_\Q$-submodules 
$$ E_{f,g}:=\ker(f^{m_f}(P_2))\cap\ker(g^{m_g}(P_4))\leq E_\Q ,$$
and their $\Q$-dimension $d_{f,g}$, where $m_f$ and $m_g$ are the associated
multiplicities in the respective characteristic polynomials.
The result is given Table \ref{eigentbl}, where (to save space) 
we abbreviate $f_4:=X^4-16X^3-75X^2+1706X-2768\in\Q[X]$ and  
$g_4:=X^4-1288X^3+405424X^2-2113152X-2701694976\in\Q[X]$.

\abs
Since we obtain $14$ pairwise distinct submodules, we conclude 
that these coincide with the homogeneous components of $E_\Q$.

\abs\bfi
The cases such that $d_{f,g}=1$ correspond to the rational constituents 
of $1_H^G$ with multiplicity $1$, where the very first case 
corresponds to the trivial character $1_G$. The action of the $P_j$
on the various $E_{f,g}$ directly yields the associated characters of $E_\C$. 
The orthogonality relations yield the degree of their Fitting
correspondents, determining the Fitting correspondence in these cases.

\abs\bfii
We consider the cases such that $d_{f,g}=2$, and compute the 
splitting fields of the irreducible divisors of degree $2$ occurring.
From this we conclude that these cases correspond
to the non-rational constituents of $1_H^G$ with multiplicity $1$:

\abs
In the first case we get simultaneous splitting field $\Q(r_3)$, 
over which $E_{f,g}$ splits into two $1$-dimensional submodules.
Thus this case corresponds to $\chi_{23/24}$,
and as above we determine the associated characters of $E_\C$. 

\abs
Similarly, the second and third cases yield simultaneous splitting field 
$\Q(r_5)$, over which both submodules split. Thus these cases 
correspond to $\chi_{36/37;38/39}$,
and as above we determine the associated characters of $E_\C$. 
The orthogonality relations yield the degree of the associated Fitting
correspondents, showing that these cases correspond to 
$\chi_{38/39}$ and $\chi_{36/37}$, respectively.

\abs\bfiii
Hence we conclude that the cases $d_{f,g}>2$ correspond to 
the remaining constituents of $1_H^G$, that is, $\chi_{19/20;21}$,
each of which occurs with multiplicity $2$. 

\abs
Since $\chi_{21}$ is rational, while $\chi_{19/20}$ is not, 
we infer that the case $d_{f,g}=4$ corresponds to $\chi_{21}$;
thus the trace map afforded by the action of the $P_j$ on $E_{f,g}$
is twice the associated character of $E_\C$. 

\abs 
Hence the case $d_{f,g}=8$ corresponds to $\chi_{19/20}$. Over $\Q(r_{33})$,
both $f_4$ and $g_4$ split into two irreducible factors of degree $2$, 
and $E_{f,g}$ splits into a direct sum of two $4$-dimensional submodules. 
The trace maps afforded by the action of the $P_j$ on the latter are 
twice the associated characters of $E_\C$.


\begin{table}
\caption{Simultaneous generalized eigenspaces of $P_2$ and $P_4$.}
\label{eigentbl}
$$ \begin{array}{|lc|lc||c|c||l|c|}
\hline 
f\spmid\text{char.pol.}(P_2) & m_f & 
g\spmid\text{char.pol.}(P_4) & m_g & 
d_{f,g} & \text{spl.} & \chi_i & m_i \\
\hline \hline
X-31 & 1 & X-17360 & 1 & 1 & &    \chi_1 & 1 \\
X-16 & 1 & X-1640  & 1 & 1 & &    \chi_8 & 1 \\
X-9  & 1 & X+20    & 1 & 1 & & \chi_{22} & 1 \\
X-5  & 1 & X+120   & 1 & 1 & & \chi_{32} & 1 \\
X-1  & 3 & X-284   & 1 & 1 & & \chi_{29} & 1 \\
     &   & X-196   & 1 & 1 & & \chi_{30} & 1 \\
     &   & X-20    & 1 & 1 & & \chi_{51} & 1 \\
X+2  & 1 & X-112   & 1 & 1 & & \chi_{14} & 1 \\
X+12 & 1 & X-1192  & 1 & 1 & & \chi_{11} & 1 \\
\hline
X^2-6X-39  & 1 & X^2+100X-1388  & 1 & 2 & r_3    & \chi_{23/24} & 1 \\
X^2-45     & 1 & X^2+130X+3820  & 1 & 2 & r_5    & \chi_{38/39} & 1 \\
X^2+12X+31 & 1 & X^2+110X-620   & 1 & 2 & r_5    & \chi_{36/37} & 1 \\
\hline
X^2+3X-64  & 2 & X^2-424X-11520 & 2 & 4 &        & \chi_{21}    & 2 \\
\hline
f_4        & 2 & g_4            & 2 & 8 & r_{33} & \chi_{19/20} & 2 \\
\hline
\end{array} $$
\absr\end{table}

\begin{table}\caption{The character table of $E_\C$.}\label{ctbl}
\tiny \setlength{\arraycolsep}{0.5em} \vspace*{185mm} \hspace*{13mm} 
\begin{rotate}{90}
$\begin{array}{|r|r||r|r||rrrrrrrrrr|} \hline
\ph&\ph^\ast&\chi_\ph&\chi_\ph(1)&1&2&3&4&5&6&8^\ast=7&7^\ast=8&9&10\\
\hline \hline
1&&1&1& 1& 31& 930& 17360& 26040& 27776& 416640& 416640& 624960& 3333120\\
2&&8&889111& 1& 16& 225& 1640& 2670& 1856& 13920& 13920& 21240& 41280\\
3&&11&1776888& 1& -12& 113& 1192& -1222& 1632& -10608& -10608& 6792& 36912\\
4&&14&4290927& 1& -2& -27& 112& 168& 672& -1008& -1008& -588& 1792\\
5&6&19&35411145& 2& 8+r_{33}& 141-10r_{33}& 644-40r_{33}& 519+31r_{33}& 
720+8r_{33}& 264+152r_{33}& 264+152r_{33}& 4422-318r_{33}& 2976-544r_{33}\\
6&5&20&35411145& 2& 8-r_{33}& 141+10r_{33}& 644+40r_{33}& 519-31r_{33}& 
720-8r_{33}& 264-152r_{33}& 264-152r_{33}& 4422+318r_{33}& 2976+544r_{33}\\
7&&21&95288172& 2& -3& 75& 424& -570& 720& 0& 0& -264& 7552\\
8&&22&230279749& 1& 9& 50& -20& 80& -120& -480& -480& -720& -320\\
9&10&23&259775040& 1& 3-4r_3& 26-24r_3& -50-36r_3& 224-8r_3& -36+88r_3& 
-564+160r_3& -564+160r_3& -108-16r_3& -872+880r_3\\
10&9&24&259775040& 1& 3+4r_3& 26+24r_3& -50+36r_3& 224+8r_3& -36-88r_3& 
-564-160r_3& -564-160r_3& -108+16r_3& -872-880r_3\\
11&&29&460559498& 1& 1& -30& 284& 0& -256& 0& 0& -864& 1152\\
12&&30&493456605& 1& 1& -30& 196& 0& -168& 0& 0& -336& 448\\
13&&32&786127419& 1& 5& -6& -120& -168& 96& 192& 192& 384& -256\\
14&15&36&885257856& 1& -6+r_5& 10-12r_5& -55-27r_5& 40+76r_5& 10+50r_5& 
90-170r_5& 90-170r_5& 110-70r_5& -640+520r_5\\
15&14&37&885257856& 1& -6-r_5& 10+12r_5& -55+27r_5& 40-76r_5& 10-50r_5&
90+170r_5& 90+170r_5& 110+70r_5& -640-520r_5\\
16&17&38/39&1016407168& 1& 3r_5& 14& -65-9r_5& -28-48r_5& 6+6r_5&
42+6r_5& 42+6r_5& -666+66r_5& 904-48r_5\\
17&16&39/38&1016407168& 1& -3r_5& 14& -65+9r_5& -28+48r_5& 6-6r_5& 
42-6r_5& 42-6r_5& -666-66r_5& 904+48r_5\\
18&&51&1842237992& 1& 1& -30& 20& 0& 8& 0& 0& 720& -960\\
\hline \end{array}$ \end{rotate}
\hspace*{45mm} \begin{rotate}{90}
$\begin{array}{|r||rrrrrrrrr|} \hline
\ph&11&13^\ast=12&12^\ast=13&14&15&16&17&19^\ast=18&18^\ast=19\\
\hline \hline
1& 4999680& 6666240& 6666240& 9999360& 13332480& 53329920& 66060288& 
79994880& 79994880\\ 
2& 77760& 49920& 49920& 72000& 96000& 61440& -73728& 230400& 230400\\
3& 21312& 28416& 28416& -23424& -50496& 29184& 12288& -137856& -137856\\
4& -19488& -1344& -1344& 23856& 66304& 146944& -172032& -45696& -45696\\
5& -3456-1328r_{33}& -6720-704r_{33}& -6720-704r_{33}& 1752+1256r_{33}& 
-384-192r_{33}& 14208+1792r_{33}& 26112+4608r_{33}&
-18816-640r_{33}& -18816-640r_{33}\\
6& -3456+1328r_{33}& -6720+704r_{33}& -6720+704r_{33}& 1752-1256r_{33}& 
-384+192r_{33}& 14208-1792r_{33}& 26112-4608r_{33}& -18816+640r_{33}&
-18816+640r_{33}\\
7& -7680& -2496& -2496& -5376& -3200& 7168& -7680& 10752& 10752\\
8& 3920& -4160& -4160& -8000& 9280& 1920& -3072& -640& -640\\
9& -2644+1024r_3& 1312-64r_3& 1312-64r_3& 736-1376r_3& 136+2224r_3& 
4704+3392r_3& 9600+8960r_3& 224+64r_3& 224+64r_3\\
10& -2644-1024r_3& 1312+64r_3& 1312+64r_3& 736+1376r_3& 136-2224r_3& 
4704-3392r_3& 9600-8960r_3& 224-64r_3& 224-64r_3\\
11& 2400& -1536& -1536& -1536& -1920& 768& 0& 3072& 3072\\
12& -2352& 2688& 2688& 2688& 5824& -6272& 0& -5376& -5376\\
13& 1728& -768& -768& 2304& 256& -6144& -6144& -6144& -6144\\
14& -690+146r_5& 320+112r_5& 320+112r_5& -2040-280r_5& 1040-8r_5& 
-3760-1440r_5& -3296-1440r_5& 3240+328r_5& 3240+328r_5\\
15& -690-146r_5& 320-112r_5& 320-112r_5& -2040+280r_5& 1040+8r_5& 
-3760+1440r_5& -3296+1440r_5& 3240-328r_5& 3240-328r_5\\
16& 878-318r_5& 592+480r_5& 592+480r_5& 1864-696r_5& -1304+528r_5& 
1776-384r_5& 2976+1632r_5& -424+2040r_5& -424+2040r_5\\
17& 878+318r_5& 592-480r_5& 592-480r_5& 1864+696r_5& -1304-528r_5&
1776+384r_5& 2976-1632r_5& -424-2040r_5& -424-2040r_5\\
18& -240& -480& -480& -480& -1920& 2880& 0& 960& 960\\
\hline \end{array}$ \end{rotate}
\hspace*{45mm} \begin{rotate}{90}
$\begin{array}{|r||rrrrrrrr|} \hline
\ph&20&21&22&23&24&25&26&27\\ \hline \hline
1& 159989760& 159989760& 319979520& 341311488& 
1279918080& 1279918080& 2047868928& 2559836160\\
2& 276480& 138240& 460800& 12288& 0& 46080& -1179648& -645120\\
3& -341760& 186624& 439296& -288768& 172032& 766464& -147456& -580608\\
4& 26880& -32256& -182784& -258048& -258048& 451584& 344064& -43008\\
5& -1248+4832r_{33}& 32928+3680r_{33}& -43392-2944r_{33}& 
86784+5376r_{33}& -60672+46336r_{33}& -108288+36608r_{33}& 
247296-19968r_{33}& -150528-76800r_{33}\\
6& -1248-4832r_{33}& 32928-3680r_{33}& -43392+2944r_{33}& 
86784-5376r_{33}& -60672-46336r_{33}& -108288-36608r_{33}& 
247296+19968r_{33}& -150528+76800r_{33}\\
7& 12480& -29376& -39168& -39936& 218112& -91392& -192000& 153600\\
8& -4800& -18880& 25600& 17920& 17920& 17920& -27648& -20480\\
9& 5856-7616r_3& 10880-3744r_3& -2048-12256r_3& 1408+4992r_3&
-32192-12288r_3& 17920+13376r_3& -78720-11520r_3& 63232+13632r_3\\
10& 5856+7616r_3& 10880+3744r_3& -2048+12256r_3& 1408-4992r_3&
-32192+12288r_3& 17920-13376r_3& -78720+11520r_3& 63232-13632r_3\\
11& 7680& 21504& -19968& -12288& -12288& 12288& 0& 0\\
12& -13440& -8064& 5376& 21504& 21504& -21504& 0& 0\\
13& -2304& 11520& -12288& 12288& 0& 12288& -30720& 30720\\
14& -4720+1696r_5& 8000+320r_5& -1320-2136r_5& -3520+256r_5& 
6160+5008r_5& 320-6336r_5& 12576+5344r_5& -15520-2208r_5\\
15& -4720-1696r_5& 8000-320r_5& -1320+2136r_5& -3520-256r_5& 
6160-5008r_5& 320+6336r_5& 12576-5344r_5& -15520+2208r_5\\
16& 2976-1200r_5& -6640+1008r_5& 1192-2760r_5& -3392-3456r_5& 
-5840-2448r_5& -6272+768r_5& 24480+8928r_5& -13280-6624r_5\\
17& 2976+1200r_5& -6640-1008r_5& 1192+2760r_5& -3392+3456r_5& 
-5840+2448r_5& -6272-768r_5& 24480-8928r_5& -13280+6624r_5\\
18& 2400& -9120& 9600& -3840& -3840& 3840& 0& 0\\
\hline \end{array}$ \end{rotate}
\end{table}

\Abs
Collecting the traces of the intersection matrices $P_1\ld P_{27}$
on the various generalized eigenspaces yields the character table 
$\Phi=[\ph_1\ld\ph_{18}]\in\C^{18\tm 27}$ of $E_\C$, see Table \ref{ctbl}.
We also indicate the pairing of $H$-orbits, Galois conjugate characters, 
the Fitting correspondence, and the degree of the Fitting correspondents.  
The rows are ordered such that the Fitting correspondents appear as in 
Table \ref{pchtbl}, except that we only know that
$\ph_{16/17}$ correspond to $\chi_{38/39}$.
 
\abs
Indeed, the character fields of the $\ph_j$ coincide with the character
fields of their respective Fitting correspondents. Since the quadratic
fields $\Q(r_3)$, $\Q(r_5)$ and $\Q(r_{33})$ are disjoint, it follows 
that there are Galois automorphisms inducing the involutions 
$(\chi_{19}\llra\chi_{20})$ and $(\chi_{23}\llra\chi_{24})$ and 
$(\chi_{36}\llra\chi_{37})(\chi_{38}\llra\chi_{39})$. 
A consideration of the ordinary character table of $G$ shows that
actually there is no table automorphism interchanging $\chi_{36/37}$
and fixing $\chi_{38/39}$. This says that by considering ordinary 
character theory alone the Fitting correspondence is only determined
up to the ambiguity above.

\abs
We define $\ph_{14/15/16,17}$ by letting 
$\ph_{14}(P_2)=-6+r_5$ and $\ph_{15}(P_2)=-6-r_5$, 
and $\ph_{16}(P_2)=3r_5$ and $\ph_{17}(P_2)=-3r_5$,
being the roots of $X^2+12X+31$ and $X^2-45$, respectively,
see Table \ref{eigentbl}. 
Then, choosing $\chi_{\ph_{14}}:=\chi_{36}$ and $\chi_{\ph_{15}}:=\chi_{37}$,
it remains to decide by further inspection whether
$\chi_{\ph_{16}}=\chi_{38}$ or $\chi_{\ph_{16}}=\chi_{39}$.

\section{Decomposition numbers}\label{sec5}

\Abs
Let now $\F:=\F_{11}$. Having the regular representation of $E_\F$ at hand, 
we apply the \MA{} to compute the simple $E_\F$-modules $S$ and their
multiplicities $[E_\F\cn S]$ in the regular $E_\F$-module. We get
$$ \begin{array}{llll} (1a)^{10},&(1b)^8,&(1c)^6,&(1d)^3, \\ \end{array} $$
saying that all constituents are $1$-dimensional, 
and where exponents denote their multiplicity. In particular, 
$\F$ already is a splitting field for $E_\F$.

\abs
Let $P_i$, where $i\in\{a,b,c,d\}$, be the projective indecomposable 
$E_\F$-module associated with the simple module $S_i:=(1i)$,
that is, such that $P_i/\rad(P_i)\cong S_i$. Since $H$ is an $11'$-group, 
we have $11\spnmid\prod_{j=1}^{27}n_j$, implying that $E_\F$ is a 
symmetric algebra, thus having a symmetric Cartan matrix and 
$\dim_\F(P_i)=[E_\F\cn S_i]$.

\abs
Using the \MA{} again to compute the indecomposable direct
summands of the regular module, we find the following Cartan matrix 
of $E_\F$, giving the multiplicities $[P_i\cn S_j]$,
both rows and columns being parameterized by $\{a,b,c,d\}$:
$$ \left[\begin{array}{rrrr} 
7 & 3 & . & . \\ 3 & 5 & . & . \\ . & . & 6 & . \\ . & . & . & 3 \\
\end{array}\right] .$$
This shows that $E_\F$ has three blocks, given by 
$\{S_a,S_b\}\dcup\{S_c\}\dcup\{S_d\}$.

\abs
Let $e_i\in E_\F$ be a primitive idempotent associated with $S_i$, 
where we may assume that the $e_i$ are pairwise orthogonal. 
Then we have $P_i\cong e_iE_\F$ as $E_\F$-modules, and 
$E_\F=(e_aE_\F\oplus e_bE_\F)\oplus e_cE_\F \oplus e_dE_\F$,
brackets indicating blocks.

\abs
Hence, applying the primitive idempotents $e_i\in E_\F$ 
to the permutation module $\F_H^G$ we get the direct sum decomposition
$\F_H^G=\F_H^G e_a\oplus\F_H^G e_b\oplus\F_H^G e_c\oplus\F_H^G e_d$
into indecomposable, pairwise non-isomorphic $\F[G]$-modules. 
In particular, this says that $P_{\F_G}$ is {\it not} a permutation module. 

\begin{table}\caption{The character table of $E_\F$.}\label{modtbl}
$$ \begin{array}{|r||rrrrrrrrrrrrrr|} \hline
(\ph_j)_\F& 1& 2& 3& 4& 5& 6& 7& 8& 9& 10& 11& 12& 13& 14 \\ \hline 
1 & 1&  9& 6& 2&  3& 1& 4& 4&  6& 10& 4& 9& 9& 8 \\
3 & 1& 10& 3& 4& 10& 4& 7& 7&  5&  7& 5& 3& 3& 6 \\
9 & 1&  1& 3& 9&  0& 8& 0& 0&  5&  8& 2& 4& 4& 4 \\
2 & 1&  5& 5& 1&  8& 8& 5& 5& 10&  8& 1& 2& 2& 5 \\ \hline 
\end{array} $$
\vspace*{0.5em}
$$ \begin{array}{|r||rrrrrrrrrrrrr|} \hline
(\ph_j)_\F& 15& 16& 17& 18& 19& 20& 21& 22& 23& 24& 25& 26& 27 \\ \hline 
1 & 7& 6& 8& 9& 9&  7&  7&  3&  1&  1& 1&  6& 2 \\
3 & 5& 1& 1& 7& 7& 10&  9&  0&  4&  3& 6& 10& 5 \\
9 & 5& 9& 0& 3& 3&  2& 10&  8& 10& 10& 1&  0& 0 \\
2 & 3& 5& 5& 5& 5&  6&  3& 10&  1&  0& 1&  3& 8 \\ \hline  
\end{array} $$ 
\absr\end{table}

\Abs
Although this already answers the question from the beginning, we can do 
better, and try and determine the projective indecomposable characters of 
$G$ being contained in $1_H^G$, in particular encompassing $\Psi_{\F_G}$.
To this end, we determine the $11$-modular decomposition matrix of $E$: 

\abs
Using $\GAP{}$, we compute the $11$-modular reduction 
$\Phi_\F:=[(\ph_1)_\F\ld(\ph_{18})_\F]\in\F^{18\tm 27}$
of the character table of $E_\C$. It turns out that 
$\{(\ph_1)_\F,(\ph_3)_\F,(\ph_9)_\F,(\ph_2)_\F\}$ are 
pairwise distinct and $\F$-linearly independent, all having degree $1$,
see Table \ref{modtbl}.
(Hence $\{\ph_1,\ph_3,\ph_9,\ph_2\}$ is a `Basic Set' 
in the sense of \cite[Def.3.1.1]{MOCBuch}.) 

\abs
Then the relations between the rows of $\Phi_\F$,
together with the fact the character degrees involved are at most $2$,
yield the complete $11$-modular decomposition matrix of $E$.
In particular, this shows that the blocks of characters consist of 
$$ \{\ph_1,\ph_3,\ph_4,\ph_5,\ph_6,\ph_7,\ph_8,\ph_{14},\ph_{17}\},
   \{\ph_9,\ph_{11},\ph_{12},\ph_{15},\ph_{16},\ph_{18}\}, 
   \{\ph_2,\ph_{10},\ph_{13}\} .$$
The blocks of the $11$-modular decomposition matrix of $E$ 
are given in Table \ref{decmat},
where we also repeat the Fitting correspondence from 
Table \ref{ctbl} and the multiplicities $m_i$ from Table \ref{pchtbl}.

\abs
Comparing the dimension of the projective indecomposable
$E_\F$-modules and the multiplicity of the simple $E_\F$-modules 
as constituents of the regular module yields the following correspondence:
$$ S_a \llra (\ph_1)_\F, \quad S_b \llra (\ph_3)_\F, \quad
   S_c \llra (\ph_9)_\F, \quad S_d \llra (\ph_2)_\F.  $$

\begin{table}\caption{The $11$-modular decomposition matrix of $E$.}
\label{decmat}
$$ \begin{array}{|r||r|r||rr|} \hline 
\ph_j & \chi_i & m_i & (\ph_1)_\F & (\ph_3)_\F \\
\hline \hline
 1 &  1    & 1 & 1 & . \\ 
 3 & 11    & 1 & . & 1 \\ 
 4 & 14    & 1 & 1 & . \\ 
 5 & 19    & 2 & 1 & 1 \\ 
 6 & 20    & 2 & 1 & 1 \\ 
 7 & 21    & 2 & 1 & 1 \\ 
 8 & 22    & 1 & 1 & . \\ 
14 & 36    & 1 & 1 & . \\ 
17 & 38/39 & 1 & . & 1 \\ 
\hline \end{array} \quad \quad \quad
\begin{array}{|r||r|r||r|} \hline 
\ph_j & \chi_i & m_i & (\ph_9)_\F \\
\hline \hline
 9 &    23 & 1 & 1 \\ 
11 &    29 & 1 & 1 \\ 
12 &    30 & 1 & 1 \\ 
15 &    37 & 1 & 1 \\ 
16 & 39/38 & 1 & 1 \\ 
18 &    51 & 1 & 1 \\ 
\hline \multicolumn{4}{c}{} \\ \hline 
\ph_j & \chi_i & m_i & (\ph_2)_\F \\
\hline \hline
 2 &     8 & 1 & 1 \\ 
10 &    24 & 1 & 1 \\ 
13 &    32 & 1 & 1 \\ 
\hline \end{array} $$
\absr\end{table}

\Abs
Thus, by Fitting correspondence, we have determined four columns
of the $11$-modular decomposition matrix of $G$, 
up to the ambiguity for the Fitting correspondents of $\ph_{16/17}$,
having fixed those of $\ph_{5/6}$, $\ph_{9/10}$ and $\ph_{14/15}$;
see Table \ref{pchagaintbl}, where $a\in\{0,1\}$.
(Actually, it is possible to decide which of the two cases left 
actually occurs, again by a computational attack similar to the
one described here; details about this will appear elsewhere \cite{part2}.)

\abs
In view of the remarks in \ref{psi}, we observe that 
the projective indecomposable character $\Psi_{\F_G}$ associated
with the trivial character, which corresponds to the projective
indecomposable $E_\F$-module $P_a$, is a non-rational character indeed.

\abs
Finally, the results of this article constitute the first steps towards 
the ambitious goal of finding the complete $11$-modular decomposition 
matrix of $G=J_4$. Richard Parker would have been keen to pursue this!

\begin{table} 
\caption{The permutation character on the cosets of $H$ in $G$ again.}
\label{pchagaintbl}
$$ \begin{array}{|r|r||rrrr|} \hline
\chi_i & \chi_i^\ast & \Psi_{\F_G} & & & \\ \hline \hline
 1 &    & 1 & . & . & . \\
 8 &    & . & 1 & . & . \\
11 &    & . & . & 1 & . \\
14 &    & 1 & . & . & . \\
19 & 20 & 1 & . & 1 & . \\
20 & 19 & 1 & . & 1 & . \\
21 &    & 1 & . & 1 & . \\
22 &    & 1 & . & . & . \\
23 & 24 & . & . & . & 1 \\
24 & 23 & . & 1 & . & . \\
29 &    & . & . & . & 1 \\
30 &    & . & . & . & 1 \\
32 &    & . & 1 & . & . \\
36 & 37 & 1 & . & . & . \\
37 & 36 & . & . & . & 1 \\
38 & 39 & . & . & a &1-a\\
39 & 38 & . & . &1-a& a \\
51 &    & . & . & . & 1 \\ 
\hline \end{array} $$
\absr\end{table}

\absr


\absr\abs\abs
{\sc
Chair for Algebra and Number Theory, RWTH Aachen University \\
Pontdriesch 14/16, D-52062 Aachen, Germany} \\
{\sf e-mail: juergen.mueller@math.rwth-aachen.de}

\end{document}